\documentclass{amsart}
\usepackage{amssymb}

\usepackage{amsmath}
\usepackage{amscd}
\usepackage{thmdefs}

\input{tcilatex}
\theoremstyle{definition}
\theoremstyle{remark}
\numberwithin{equation}{section}

\input tcilatex

\begin{document}
\title[Moments of Operators in $L(F_{N})$]{The Moment of an Operator in the Free Group Factor $L(F_{N})$}
\author{Ilwoo Cho}
\address{Dep. of Math, Univ. of Iowa, Iowa City, IA, U. S. A}
\email{ilcho@math.uiowa.edu}
\maketitle

\begin{abstract}
In this paper, we will define an operator $X_{n}$ by the total sum of all
word with their length $n$ such that $X_{n}=\underset{\left| w\right| =n}{%
\sum }w$ in the free group factor $L(F_{N}),$ where $F_{N}$ is the free
group with $N$-generators. We will construct the recurrence relation of the
operator product $x^{k}X_{n},$ where $x=X_{1}$ is the generating operator of 
$L(F_{N}),$ for $k,$ $n$ $\in $ $\Bbb{N}.$ By this recurrence relation, we
can compute the moment $\tau \left( x^{k}X_{n}\right) $ of $x^{k}X_{n},$ for
the cases when $k=n$ and $k<n,$ where $\tau $ $:$ $L(F_{N})$ $\rightarrow $ $%
\Bbb{C}$ is the canonical trace on the free group factor $L(F_{N}).$
\end{abstract}

\strut

From mid 1980's, Free Probability Theory has been developed. Here, the
classical concept of Independence in Probability theory is replaced by a
noncommutative analogue called Freeness (See [9]). There are two approaches
to study Free Probability Theory. One of them is the original analytic
approach of Voiculescu and the other one is the combinatorial approach of
Speicher and Nica (See [1], [2] and [3]). \medskip Let $A$ be a von Neumann
algebra and let $\varphi :A\rightarrow \Bbb{C}$ be a \ linear functional
satisfying that $\varphi (a^{*})$ $=$ $\overline{\varphi (a)},$ for all $a$ $%
\in $ $A.$ Then the algebraic pair $\left( A,\varphi \right) $ is called the 
$W^{*}$-probability space. All elements in $\left( A,\varphi \right) $ are
called random variables. The basic free probabilistic information of the
fixed random variable $a$ $\in $ $\left( A,\varphi \right) $ is the (free)
moments $\varphi (a^{n}),$ for $n$ $\in $ $\Bbb{N},$ of the random variable $%
a$. Throughout this paper, let

$\strut $

\begin{center}
$F_{N}=\,<g_{1},$ $g_{2},$ $...,$ $g_{N}>$
\end{center}

\strut

be the free group with $N$-generators. Then we can construct the free group
factor $L(F_{N}).$ i.e.,

\strut

\begin{center}
$L(F_{N})=\overline{\Bbb{C}[F_{N}]}^{w}.$
\end{center}

\strut

This von Neumann algebra is indeed a factor, because the group $F_{N}$ is
icc. (Recall that the group von Neumann algebra $L(G)$ is a factor if and
only if the group $G$ is an icc group.) Let $a$ be an operator in $L(F_{N}).$
Then there exists the Fourier expansion of $x,$

\strut

\begin{center}
$a=\underset{g\in F_{N}}{\sum }\alpha _{g}u_{g},$ \ with $\alpha _{g}\in %
\mathbb{C},$ for all $g\in F_{N}.$
\end{center}

\strut

We can regard all $g\in F_{N}$ as unitaries $u_{g}$ in $L(F_{N}).$ For the
convenience, we will denote these unitaries $u_{g}$ just by $g.$ With this
notation, it is easy to check that

$\strut $

\begin{center}
$g^{*}=u_{g}^{*}=u_{g}^{-1}=u_{g^{-1}}=g^{-1}$ in $L(F_{N}),$
\end{center}

\strut

where $g^{-1}$ is the group inverse of $g$ in $F_{N}.$ We can define the
canonical trace $\tau :L(F_{N})\rightarrow \Bbb{C}$ by

\strut

\begin{center}
$\tau \left( \underset{g\in F_{2}}{\sum }\alpha _{g}g\right) =\alpha _{e}.$
\end{center}

\strut

Then the algebraic pair $\left( L(F_{N}),\tau \right) $ is a $W^{*}$%
-probability space. \strut

\strut \strut

In [15], we re-computed the moments of the generating operator

\strut

\begin{center}
$x=\sum_{j=1}^{N}\left( g_{j}+g_{j}^{-1}\right) $
\end{center}

\strut

of $L(F_{N}),$ by using the following recurrence diagram,

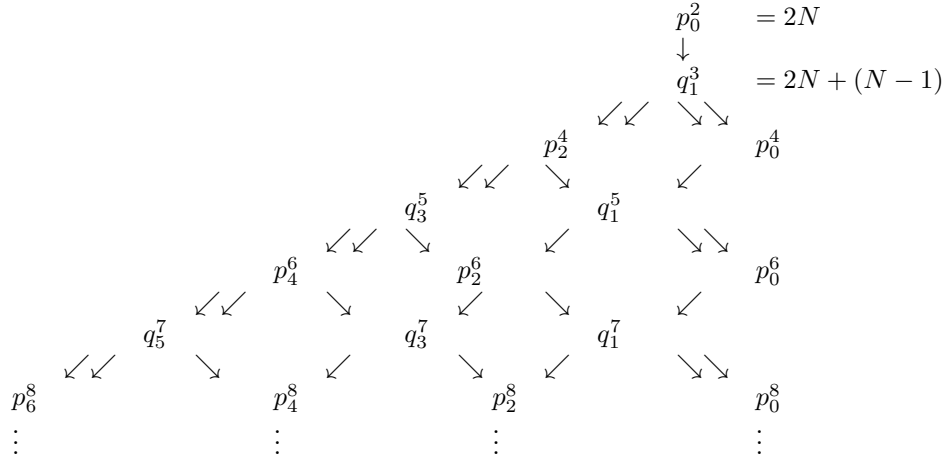
\strut

\begin{center}
$
\begin{array}{llllllllllll}
&  &  &  &  &  &  &  &  &  & p_{0}^{2} & =2N \\ 
&  &  &  &  &  &  &  &  &  & \downarrow &  \\ 
&  &  &  &  &  &  &  &  &  & q_{1}^{3} & =2N+(N-1) \\ 
&  &  &  &  &  &  &  &  & \swarrow \swarrow & \searrow \searrow &  \\ 
&  &  &  &  &  &  &  & p_{2}^{4} &  &  & p_{0}^{4} \\ 
&  &  &  &  &  &  & \swarrow \swarrow & \searrow &  & \swarrow &  \\ 
&  &  &  &  &  & q_{3}^{5} &  &  & q_{1}^{5} &  &  \\ 
&  &  &  &  & \swarrow \swarrow & \searrow &  & \swarrow &  & \searrow
\searrow &  \\ 
&  &  &  & p_{4}^{6} &  &  & p_{2}^{6} &  &  &  & p_{0}^{6} \\ 
&  &  & \swarrow \swarrow &  & \searrow &  & \swarrow & \searrow &  & 
\swarrow &  \\ 
&  & q_{5}^{7} &  &  &  & q_{3}^{7} &  &  & q_{1}^{7} &  &  \\ 
& \swarrow \swarrow &  & \searrow &  & \swarrow &  & \searrow & \swarrow & 
& \searrow \searrow &  \\ 
p_{6}^{8} &  &  &  & p_{4}^{8} &  &  & \text{ \ \ \ }p_{2}^{8} &  &  &  & 
p_{0}^{8} \\ 
\vdots &  &  &  & \vdots &  &  & \text{ \ \ \ }\vdots &  &  &  & \vdots
\end{array}
$
\end{center}

\strut

where

\begin{center}
$\swarrow \swarrow $ \ : \ $(2N-1)+[$former term$]$

$\searrow $ \ \ \ \ : \ \ $(2N-1)\cdot [$former term$]$

$\swarrow $ \ \ \ \ : \ \ $\cdot +[$former term$]$
\end{center}

and

\begin{center}
$\searrow \searrow $ \ : \ \ $(2N)\cdot [$former term$].$
\end{center}

\strut \strut

The numbers in the recurrence diagram are came from the well-known relations

\strut

\begin{center}
$X_{1}X_{1}=X_{2}+2Ne$
\end{center}

and

\begin{center}
$X_{1}X_{n}=X_{n+1}+(N-1)X_{n-1},$
\end{center}

\strut

for all $N,\,\,n\in \Bbb{N}\,\setminus \,\{1\}$ (See [16]). For example,
since $x=X_{1}$ and $x^{3}$ $=$ $X_{1}$ $X_{1}$ $X_{1},$ we can have

\strut

\begin{center}
$
\begin{array}{ll}
x^{3} & =X_{1}\left( X_{2}+2Ne\right) =X_{1}X_{2}+2N\cdot X_{1} \\ 
&  \\ 
& =X_{3}+\left( (N-1)+2N\right) X_{1}=X_{3}+q_{1}^{3}X_{1}.
\end{array}
$
\end{center}

\strut

This recurrence diagram represents that

\strut

\begin{center}
$x^{2k}=X_{2k}+p_{2k-2}^{2k}X_{2k-2}+...+p_{2}^{2k}X_{2}+p_{0}^{2k}e$
\end{center}

and

\begin{center}
$%
x^{2k+1}=X_{2k+1}+q_{2k-1}^{2k+1}X_{2k-1}+...+q_{3}^{2k}X_{3}+q_{1}^{2k}X_{1}, 
$
\end{center}

\strut

for all $k\in \Bbb{N},$ where $p_{j}^{2k}$'s and $q_{i}^{2k+1}$'s are gotten
from the above recurrence diagram and where $e$ is the identity of $F_{N}$
and where $X_{n}$ $=$ $\underset{\left| w\right| =n}{\sum }w$ is the total
sum of all words with their length $n,$ as an operator in the free group
factor $L(F_{N}),$ for all $n\in \Bbb{N}.$ Therefore, we can get that all
odd moments of the generating operator $x$ of $L(F_{N})$ vanish. More
precisely, we have that

\strut

\begin{center}
\strut $\tau \left( x^{n}\right) =\left\{ 
\begin{array}{lll}
0 &  & \text{if }n\text{ is odd} \\ 
&  &  \\ 
p_{0}^{n} &  & \text{if }n\text{ is even,}
\end{array}
\right. $
\end{center}

\strut

where $p_{0}^{n}$'s are gotten from the above recurrence diagram, for all $%
n\in 2\Bbb{N}.$

\strut

In this paper, we will consider the operators $x^{k}X_{n},$ in $L(F_{N}),$
for $k,$ $n$ $\in $ $\Bbb{N}.$ By regarding them as random variables in the $%
W^{*}$-probability space $\left( L(F_{N}),\tau \right) ,$ we can compute the
moment $\tau (x^{k}X_{n})$. In order to do that, we will construct another
recurrence relation to express $x^{k}X_{n},$ in terms of $X_{j}$'s. This
recurrence relation is needed because there is no concrete recurrence
relations for $X_{m}X_{n},$ where $m,$ $n$ $\in $ $\Bbb{N}$ $\setminus $ $%
\{1\}.$ One of the main results of this paper is that if $k=n,$ then 

\strut 

\begin{center}
$
\begin{array}{ll}
\tau \left( x^{n}X_{n}\right)  & =\tau \left(
r_{1}^{(n)}X_{2n}+r_{2}^{(n)}X_{2n-2}+...+r_{n}^{(n)}X_{2}+r_{n+1}^{(n)}e%
\right)  \\ 
&  \\ 
& =r_{n+1}^{(n)},
\end{array}
$
\end{center}

\strut 

where the sequence $\left( r_{1}^{(n)},...,r_{n+1}^{(n)}\right) $ is the
coefficient sequence of $\left( r+(N-1)\right) ^{n},$ for all $n\in \Bbb{N}.$
Here, $r$ is just an indeterminant. Also, it is shown that if $n>k,$ then $%
\tau (x^{k}X_{n})=0.$

\strut \strut

\strut \strut

\strut \strut

\section{\strut The Operator $x^{k}X_{n}$ in $L(F_{N})$}

\strut

\strut

Let's consider the coefficient of $\left( r+(N-1)\right) ^{n},$ for $n\in 
\Bbb{N},$ where $r$ is an arbitrary indeterminant. Then we have that the
Pascal's triangle expressing the coefficients of $\left( r+(N-1)\right) ^{n},
$ as follows ;

\strut $\strut $

\begin{center}
$
\begin{array}{lllllllll}
&  &  & \,\,\,\,\,\,\,\,1 &  &  &  & \longrightarrow  & \left(
r+(N-1)\right) ^{0} \\ 
&  & \,\,\,\,\,\,\,\,1 &  & \,\,N-1 &  &  & \longrightarrow  & \left(
r+(N-1)\right) ^{1} \\ 
& 1 &  & 2(N-1) &  & (N-1)^{2} &  & \longrightarrow  & \left( r+(N-1)\right)
^{2} \\ 
1 &  & 3(N-1) &  & 3(N-1)^{2} &  & (N-1)^{3} & \longrightarrow  & \left(
r+(N-1)\right) ^{3} \\ 
&  &  &  &  &  &  &  &  \\ 
\vdots  &  &  &  &  &  & \,\,\,\,\,\,\,\,\,\,\vdots  &  & 
\,\,\,\,\,\,\,\,\,\,\,\,\,\,\,\,\,\,\,\,\,\vdots 
\end{array}
$
\end{center}

\strut \strut 

\strut 

\begin{definition}
We will denote the coefficients of $\left( r+(N-1)\right) ^{n}$ by the
sequence $(r_{1}^{(n)},...,r_{n+1}^{(n)}),$ for all $n\in \Bbb{N}.$ The
sequence $(r_{1}^{(n)},...,r_{n+1}^{(n)})$ is called the coefficient
sequence of $\left( r+(N-1)\right) ^{n},$ for each $n\in \Bbb{N}.$ Remark
that, in all cases, $r_{1}^{(n)}=1.$
\end{definition}

\strut

For example, the coefficient sequence of $\left( r+(N-1)\right) ^{3}$ is

\strut

\begin{center}
$\left( 1,\,3(N-1),\,3(N-1)^{2},\,(N-1)^{3}\right) .$
\end{center}

\strut \strut

In this section, we will find the recurrence relation for $x^{k}X_{n},$
where $k,n\in \Bbb{N}.$ Observe that, since $x=X_{1},$ in our case, we have
that

\strut

$\ x^{k}X_{n}=x^{k-1}xX_{n}=x^{k-1}X_{1}X_{n}$

\strut

$\ \ =x^{k-1}\left( X_{n+1}+(N-1)X_{n-1}\right) $

\strut

$\ \ =x^{k-2}\left( X_{1}(X_{n+1}+(N-1)X_{n-1})\right) $

\strut

$\ \ =x^{k-2}\left( X_{1}X_{n+1}+(N-1)X_{1}X_{n-1}\right) $

\strut

$\ \ =x^{k-2}\left( (X_{n+2}+(N-1)X_{n})+(N-1)(X_{n}+(N-1)X_{n-2})\right) $

\strut

$\ \ =x^{k-2}\left( X_{n+2}+2(N-1)X_{n}+(N-1)^{2}X_{n-2}\right) $

\strut

$\ \ =x^{k-3}\left(
X_{1}X_{n+2}+2(N-1)X_{1}X_{n}+(N-1)^{2}X_{1}X_{n-2}\right) $

\strut

$\ \ =x^{k-3}\left(
X_{n+3}+3(N-1)X_{n+1}+3(N-1)^{2}X_{n-1}+(N-1)^{3}X_{n-3}\right) $

\strut

$\ \ =....$

\strut

Thus we can easily verify that ;

\strut

\begin{theorem}
Let $k$ and $p$ be numbers in $\Bbb{N}$ satisfying that $p\leq k.$ Then $%
x^{k}X_{n}$ is

\strut 

(1.1) $\ \ \ x^{k-p}\left(
r_{1}^{(p)}X_{n+p}+r_{2}^{(p)}X_{n+p-2}+...+r_{p}^{(p)}X_{n-p+2}+r_{p+1}^{(p)}X_{n-p}\right) ,
$

\strut 

where $(r_{1}^{(p)},...,r_{p+1}^{(p)})$ is the coefficient sequence of $%
\left( r+(N-1)\right) ^{p}.$ $\square $
\end{theorem}

\strut

The above theorem is proved by the induction on $p,$ after taking the
sufficiently big $k$. By the previous theorem, we have that ;

\strut

\begin{corollary}
(1) If $k<n$ in $\Bbb{N},$ then we have that

\strut 

$\ \ \ \
x^{k}X_{n}=r_{1}^{(k)}X_{n+k}+r_{2}^{(k)}X_{n+k-2}+...+r_{k}^{(k)}X_{n-k+2}+r_{k+1}^{(k)}X_{n-k}.
$

\strut 

(2) If $k>n$ in $\Bbb{N},$ then we have that

\strut 

$\ \ \ \ x^{k}X_{n}=x^{k-n}\left(
r_{1}^{(n)}X_{2n}+r_{2}^{(n)}X_{2n-2}+...+r_{n}^{(n)}X_{2}+r_{n+1}^{(n)}e%
\right) .$

\strut 

(3) If $k=n$ in $\Bbb{N},$ then we have that

\strut 

$\ \ \ \ \ \ \ \
x^{n}X_{n}=r_{1}^{(n)}X_{2n}+r_{2}^{(n)}X_{2n-2}+...+r_{n}^{(n)}X_{2}+r_{n+1}^{(n)}e.
$
\end{corollary}

\strut

\begin{proof}
It is easy to prove (3), by (1.1). Now, assume that $k=n+k^{\prime },$ for
some $k^{\prime }\in \Bbb{N}.$ Then, by (3), we can verify the result of
(2). Similarly, we can get (1).
\end{proof}

\strut

\begin{example}
Let $k=3$ and $n=3.$ Then, by the previous corollary, we have that

\strut 

$\ \ \ \ \ \ \
x^{3}X_{3}=r_{1}^{(3)}X_{6}+r_{2}^{(3)}X_{4}+r_{3}^{(3)}X_{2}+r_{4}^{(3)}e,$

\strut 

where $\left( r_{1}^{(3)},r_{2}^{(3)},r_{3}^{(3)},r_{4}^{(3)}\right) $ is
the coefficient sequence of $\left( r+(N-1)\right) ^{3}.$ Now, take $k=2$
and $n=3.$ Then

\strut 

$\ \ \ \ \ \ 
\begin{array}{ll}
x^{2}X_{3} & =x\left( X_{1}X_{3}\right) =X_{1}\left( X_{4}+(N-1)X_{2}\right) 
\\ 
&  \\ 
& =X_{5}+2(N-1)X_{3}+(N-1)^{2}X_{1}.
\end{array}
$

\strut 

Now, we will take $k=5$ and $n=3.$ Then

\strut 

$\ \ \ \ \ \ 
\begin{array}{ll}
x^{5}X_{3} & =x^{2}\left( x^{3}X_{3}\right)  \\ 
&  \\ 
& =x^{2}\left(
r_{1}^{(3)}X_{6}+r_{2}^{(3)}X_{4}+r_{3}^{(3)}X_{2}+r_{4}^{(3)}e\right) .
\end{array}
$

\strut 

Notice that, in the above formula, we can keep ding our process as follows ;

\strut 

$\ x^{2}\left(
r_{1}^{(3)}X_{6}+r_{2}^{(3)}X_{4}+r_{3}^{(3)}X_{2}+r_{4}^{(3)}e\right) $

\strut 

$\ \ \
=r_{1}^{(3)}x^{2}X_{6}+r_{2}^{(3)}x^{2}X_{4}+r_{3}^{(3)}x^{2}X_{2}+r_{4}^{(3)}x^{2}
$

\strut 

$\ \ \ =r_{1}^{(3)}\left(
r_{1}^{(2)}X_{8}+r_{2}^{(2)}X_{6}+r_{3}^{(2)}X_{4}\right) $

\strut 

$\ \ \ \ \ \ \ \ \ \ \ +r_{2}^{(3)}\left(
r_{1}^{(2)}X_{6}+r_{2}^{(2)}X_{4}+r_{3}^{(2)}X_{2}\right) $

\strut 

$\ \ \ \ \ \ \ \ \ \ \ +r_{3}^{(3)}\left(
r_{1}^{(2)}X_{4}+r_{2}^{(2)}X_{2}+r_{3}^{(2)}e\right) $

\strut 

$\ \ \ \ \ \ \ \ \ \ \ +r_{4}^{(3)}\left( X_{2}+(2N)e\right) $

\strut 

$\ \ \ =\left( r_{1}^{(3)}r_{1}^{(2)}\right) X_{8}+\left(
r_{1}^{(3)}r_{2}^{(2)}+r_{2}^{(3)}r_{1}^{(2)}\right) X_{6}$

\strut 

$\ \ \ \ \ \ \ \ \ \ \ +\left(
r_{1}^{(3)}r_{3}^{(2)}+r_{2}^{(3)}r_{2}^{(2)}+r_{3}^{(3)}r_{1}^{(2)}\right)
X_{4}$

\strut 

$\ \ \ \ \ \ \ \ \ \ \ +\left( \allowbreak
r_{2}^{(3)}r_{3}^{(2)}+r_{3}^{(3)}r_{2}^{(2)}+r_{4}^{(3)}\right) X_{2}$

\strut 

$\ \ \ \ \ \ \ \ \ \ \ +\left( r_{3}^{(3)}r_{3}^{(2)}+\allowbreak \left(
2N\right) r_{4}^{(3)}\right) e.$
\end{example}

\strut

\strut

\strut

\strut

\section{The Moment of $x^{k}X_{n}$}

\strut

\strut

\strut

In this chapter, we will compute the moments of the random variable $%
x^{k}X_{n}$ in our $W^{*}$-probability space $\left( L(F_{N}),\tau \right) .$
Remark that to compute the tracial value $\tau (a)$ of an arbitrary random
variable $a$ in $L(F_{N})$ is to find coefficient of $e$-term of $a.$ So, we
will try to find the $e$-term of operator $x^{k}X_{n}.$ 

\strut \strut 

\begin{theorem}
Let $k,n\in \Bbb{N}$ and let $x=X_{1}$ be the generating operator of the
free group factor $L(F_{N}).$ If $X_{n}=\underset{\left| w\right| =n}{\sum }w
$ in $L(F_{n}),$ then

\strut 

(1) $\ \ \tau (x^{k}X_{n})=r_{n+1}^{(n)},$ \ \ \ whenever $n=k.$\strut 

(2) \ $\tau (x^{k}X_{n})=0,$ \ \ \ \ \ \ \ \ whenever $n>k.$
\end{theorem}

\strut

\begin{proof}
Assume that $k=n.$ Then, by (3) of the previous corollary, we have that 

\strut 

$\ \ \ \ \ \ 
\begin{array}{ll}
\tau \left( x^{n}X_{n}\right)  & =\tau \left(
r_{1}^{(n)}X_{2n}+r_{2}^{(n)}X_{2n-2}+...+r_{n}^{(n)}X_{2}+r_{n+1}^{(n)}e%
\right)  \\ 
&  \\ 
& =r_{n+1}^{(n)},
\end{array}
$

\strut 

where $\left( r_{1}^{(n)},...,r_{n+1}^{(n)}\right) $ is the coefficient
sequence of $\left( r+(N-1)\right) ^{n}.$ Now, assume that $n>k.$ Then, by
(1) of the previous corollary, we can get that

\strut

$\ 
\begin{array}{ll}
\tau \left( x^{k}X_{n}\right) & =\tau \left(
r_{1}^{(k)}X_{n+k}+r_{2}^{(k)}X_{n+k-2}+...+r_{k}^{(k)}X_{n-k+2}+r_{k+1}^{(k)}X_{n-k}\right)
\\ 
&  \\ 
& =0,
\end{array}
$

\strut \strut

since $x^{k}X_{n}$ does not have the $e$-term.
\end{proof}

\strut 

Now, suppose that $k>n$ and $k=n+k^{\prime }.$ Then

\strut

(3.1)

\begin{center}
$
\begin{array}{ll}
x^{k}X_{n} & =x^{k^{\prime }}x^{n}X_{n} \\ 
&  \\ 
& =x^{k^{\prime }}\left(
r_{1}^{(n)}X_{2n}+r_{2}^{(n)}X_{2n-2}+...+r_{n}^{(n)}X_{2}+r_{n+1}^{(n)}e%
\right) .
\end{array}
$
\end{center}

\strut

Suppose that $k^{\prime }=1.$ Then $x^{k^{\prime }}=x=X_{1}.$ So, $%
x^{k^{\prime }}X_{j}$ $=$ $X_{1}X_{j}$ does not contain $e$-term, for each $%
j $ $=$ $0,$ $2,$ $4,$ $...,$ $2n.$ This shows that

\strut

(3.2) \ \ \ \ \ \ \ \ \ \ \ if $k^{\prime }=1,$ then $\tau \left(
x^{k^{\prime }}X_{n}\right) =0.$

\strut

By (2) of the previous theorem, we have that

\strut 

\strut (3.3)

\begin{center}
$
\begin{array}{ll}
\tau \left( x^{k^{\prime }}x^{n}X_{n}\right) & =r_{1}^{(n)}\tau \left(
x^{k^{\prime }}X_{2n}\right) +r_{2}^{(n)}\tau \left( x^{k^{\prime
}}X_{2n-2}\right) \\ 
&  \\ 
& \,\,\,\,\,+...+r_{n}^{(n)}\tau \left( x^{k^{\prime }}X_{2}\right)
+r_{n+1}^{(n)}\tau \left( x^{k^{\prime }}\right) .
\end{array}
$
\end{center}

\strut

Let $j\in \{2,4,...,2n-2,2n\}$ and assume that $k^{\prime }<j.$ Then the
summands in (3.3) satisfy that

\strut

(3.4) $\ \ \ \tau \left( x^{k^{\prime }}X_{j}\right) =\tau \left(
x^{k^{\prime }}X_{j+2}\right) =...=\tau \left( x^{k^{\prime }}X_{2n}\right)
=0,$

\strut

by (2) of the previous theorem. So, we can conclude that ;

\strut

\begin{proposition}
Let $k,k^{\prime },n\in \Bbb{N}$ and $k=n+k^{\prime }$ and let $k^{\prime
}<2n.$ Assume that $j$ is the minimal number satisfying that $k^{\prime }<j,$
where $j$ $\in $ $\{2,$ $4,$ $...,$ $2n-2,$ $2n\}.$ Then

\strut 

$\ \ \ \ \ \ \ \ 
\begin{array}{ll}
\tau \left( x^{k}X_{n}\right)  & =r_{n_{j}}^{(n)}\tau \left( x^{k^{\prime
}}X_{j}\right) + \\ 
&  \\ 
& \,\,\,\,\,\,\,...+r_{n}^{(n)}\tau \left( x^{k^{\prime }}X_{2}\right)
+r_{n+1}^{(n)}\tau \left( x^{k^{\prime }}\right) .
\end{array}
$
\end{proposition}

$\strut $

\begin{proof}
Since $k>n,$ we have that

\strut

$\ \tau \left( x^{k}X_{n}\right) =\tau \left( x^{k^{\prime
}}x^{n}X_{n}\right) $

\strut

$\ \ \ =\tau \left( x^{k^{\prime
}}(r_{1}^{(n)}X_{2n}+r_{2}^{(n)}X_{2n-2}+...+r_{n}^{(n)}X_{2}+r_{n+1}^{(n)}e)\right) 
$

\strut

$\ \ \ =r_{1}^{(n)}\tau \left( x^{k^{\prime }}X_{2n}\right) +r_{2}^{(n)}\tau
\left( x^{k^{\prime }}X_{2n-2}\right) +$

\strut

$\ \ \ \ \ \ \ \ \ \ \ \ \ \ \ \ \ \ \ \ \ \ \ \ \ \ ...+r_{n}^{(n)}\tau
\left( x^{k^{\prime }}X_{2}\right) +r_{n+1}^{(n)}\tau \left( x^{k^{\prime
}}\right) .$

\strut

Assume that $j$ is the minimal number satisfying $k^{\prime }<j,$ where $%
j\in \{2,4,...,2n\}.$ Then, by (3.4), we can get that

\strut

(3.5) \ 

\ $\ \ \ \ \ \ \ 
\begin{array}{ll}
\tau \left( x^{k}X_{n}\right) & =r_{n_{j}}^{(n)}\tau \left( x^{k^{\prime
}}X_{j}\right) + \\ 
&  \\ 
& \,\,\,\,\,\,\,...+r_{n}^{(n)}\tau \left( x^{k^{\prime }}X_{2}\right)
+r_{n+1}^{(n)}\tau \left( x^{k^{\prime }}\right) .
\end{array}
$

\strut \strut
\end{proof}

\strut

Note that in the previous proposition, $\tau \left( x^{k^{\prime }}\right) $
can be computed by the recurrence diagram introduced at the beginning of
this paper. The case when $k^{\prime }>2n$ is very hard to find the concrete
formula. However, we can verify that we might have the recursive algorithm
for the computation. i.e.,

\strut

\begin{center}
$k=n+k^{\prime }=n+\left( 2n+k^{\prime \prime }\right) .$
\end{center}

\strut

So, like the observation for $k^{\prime },$ we can do the similar process
for $k^{\prime \prime }.$ Also, if $k^{\prime \prime }>4n,$ then

\strut

\begin{center}
$k=n+\left( 2n+(4n+k^{\prime \prime \prime })\right) .$
\end{center}

\strut

So, we do the similar job for $k^{\prime \prime \prime }.$ \strut 

\strut 

Notice that, for $k,n\in \Bbb{N}.$ we have that 

$\strut $

\begin{center}
$x^{k}X_{n}=X_{n}x^{k}.$
\end{center}

\strut 

Recall that $X_{1}X_{n}=X_{n}X_{1},$ for all $n\in \Bbb{N}$ in $L(F_{N}).$
By the definition of the generating operator $x$, $x=X_{1}.$ So, we have $%
xX_{n}=X_{n}x.$ Hence,

\strut 

\begin{center}
$x^{k}X_{n}=x^{k=1}xX_{n}=x^{k-1}X_{n}x=...=xX_{n}x^{k-1}=X_{n}x^{k}.$
\end{center}

\strut 

Therefore, we can get that ;

\strut 

\begin{proposition}
$\tau \left( x^{k}X_{n}\right) =\tau \left( X_{n}x^{k}\right) ,$ for all $%
k,n\in \Bbb{N}.$ $\square $
\end{proposition}

\strut 

\begin{corollary}
$\tau \left( x^{k_{1}}X_{n}x^{k_{2}}\right) =\tau \left(
x^{k_{1}+k_{2}}X_{n}\right) ,$ for all $k_{1},k_{2},n\in \Bbb{N}.$ $\square $
\end{corollary}

\strut 

\begin{corollary}
$\tau \left(
x^{k_{1}}X_{n_{1}}x^{k_{2}}X_{n_{2}}...x^{k_{m}}X_{n_{m}}\right) =\tau
\left( x^{\sum_{i=1}^{m}k_{i}}\cdot \Pi _{j=1}^{m}X_{n_{j}}\right) .$ \ $%
\square $
\end{corollary}

\strut 

Actually the above computation would be very complicated because we do not
know the concrete expression for $X_{m}X_{n},$ for all $m,$ $n$ $\in $ $\Bbb{%
N}\,\,\setminus \,\{1\}.$

\strut \strut 

\strut \strut 

\strut \strut 

\strut 

\strut \textbf{References}

\bigskip

\strut

{\small [1] \ R. Speicher, Combinatorial Theory of the Free Product with
Amalgamation and Operator-Valued Free Probability Theory, AMS Mem, Vol 132 ,
Num 627 , (1998).}

{\small [2] \ \ A. Nica, R-transform in Free Probability, IHP course note.}

{\small [3] \ \ R. Speicher, Combinatorics of Free Probability Theory IHP
course note.}

{\small [4] \ \ A. Nica, D. Shlyakhtenko and R. Speicher, R-cyclic Families
of Matrices in Free Probability, J. of Funct Anal, 188 (2002),
227-271.\strut }

{\small [5] \ \ A. Nica and R. Speicher, R-diagonal Pair-A Common Approach
to Haar Unitaries and Circular Elements, (1995), Preprint.}

{\small [6] \ \ D. Shlyakhtenko, Some Applications of Freeness with
Amalgamation, J. Reine Angew. Math, 500 (1998), 191-212.\strut }

{\small [7] \ \ A. Nica, D. Shlyakhtenko and R. Speicher, R-diagonal
Elements and Freeness with Amalgamation, Canad. J. Math. Vol 53, Num 2,
(2001) 355-381.\strut }

{\small [8] \ \ A. Nica, R-transforms of Free Joint Distributions and
Non-crossing Partitions, J. of Func. Anal, 135 (1996), 271-296.\strut }

{\small [9] \ \ D.Voiculescu, K. Dykemma and A. Nica, Free Random Variables,
CRM Monograph Series Vol 1 (1992).\strut }

{\small [10] D. Voiculescu, Operations on Certain Non-commuting
Operator-Valued Random Variables, Ast\'{e}risque, 232 (1995), 243-275.\strut 
}

{\small [11] D. Shlyakhtenko, A-Valued Semicircular Systems, J. of Funct
Anal, 166 (1999), 1-47.\strut }

{\small [12] I. Cho, The Moment Series and The R-transform of the Generating
Operator of }$L(F_{N}),${\small \ (2003), Preprint.\strut }

{\small [13] I. Cho, The Moment Series of the Generating Operator of }$%
L(F_{2})*_{L(K)}L(F_{2})${\small , (2003), Preprint.}

{\small [14] I. Cho, An Example of Moment Series under the Compatibility,
(2003), Preprint. }

{\small [15] I. Cho, Operator-Valued Moment Series of the Generating
Operator of }$L(F_{2})${\small \ over the Commutator Group von Neumann
Algebra }$L(K),$ {\small (2004), Preprint.}

{\small [16] F. Radulescu, Singularity of the Radial Subalgebra of }$%
L(F_{N}) ${\small \ and the Puk\'{a}nszky Invariant, Pacific J. of Math,
vol. 151, No 2 (1991)\strut , 297-306.\strut \strut }

\strut

\end{document}